\documentclass[12pt]{amsart}
\usepackage{amsmath,amsthm,amssymb,latexsym}
\usepackage[applemac]{inputenc}
\usepackage{wasysym}
\usepackage{phonetic}
\usepackage{palatino}
\usepackage[dvips]{color}
\usepackage[dvips]{graphicx}
\usepackage{pb-diagram}
\usepackage{pst-all}
\usepackage{pstcol}
\usepackage{tipa}
\usepackage[T1]{fontenc}
\usepackage{MnSymbol}

\newcommand\M[1]{\mathfrak{#1}}

\newcommand\C[1]{\mathcal{#1}}

\newtheorem{proposition}{Proposition}[section]
\newtheorem{theorem}[proposition]{Theorem}

\newtheorem{corollary}[proposition]{Corollary}
\newtheorem{fact}[proposition]{Fact}
\newtheorem{claim}[proposition]{Claim}
\theoremstyle{definition}
\newtheorem{definition}[proposition]{Definition}

\newtheorem{examples}[proposition]{Examples}
\newtheorem{remark}[proposition]{Remark}

\def\E{\varepsilon}

\def\ordencea{\prec_{\C{K}}}

\newcommand{\eop}[1]{
\vspace{-5.4mm}
\begin{flushright}
\qedsymbol$_{\text{#1}}$\\
\end{flushright}
}
\newenvironment{prueba}[1][{\it Proof}]{\noindent {\it #1.}}{}
\newcommand{\bdem}[1][Proof]{\begin{prueba}[#1]}
\newcommand{\edem}[1][]{\\ \eop{#1}
\end{prueba}}
\def\bsdem{\begin{prueba}[Reference]} 

\begin{document}
\title[Towards a proof of the Shelah Presentation Theorem $\dots$]{Towards a proof of the
Shelah Presentation Theorem in Metric Abstract Elementary Classes.}
\author{Pedro Zambrano}
\address{{\rm E-mail:} {\it phzambranor@unal.edu.co, phzambranor@gmail.com}\\
Departamento de Matem\'aticas Universidad Nacional de Colombia, AK 30
$\#$ 45-03 111321, Bogot\'a - Colombia}
\date{\today}
\thanks{\emph{AMS Subject Classification}: Primary: 03C52, 03C75, 03C30.
Secondary: 03C05, 03C65
  and 03C95.\\}

\begin{abstract}
In \cite{Za}, we proved the following version of Shelah's Presentation Theorem in the setting of Metric Abstract Elementary Classes:
\\ \\
{\bf Theorem.} Given an MAEC $(\C{K},\ordencea)$, there exist an
expansion $L'$ of $L(\C{K})$, an $L'$-theory $T'$ and a set of
$T'$-types $\Gamma$ (in the setting of Continuous Logic) such that $\mathcal{K}=PC_L(T',\Gamma)$ (i.e.:
$\C{K}$ is a projective class with omitting types).
\\ \\
In \cite{Za}, we claimed that the new function symbols are not necessarily uniformly continuous.
In this paper we provide a proof they are in fact uniformly continuous.\end{abstract}

\maketitle
% Eliminado el 21 de enero de 2008 ... es irrelevante, pues no voy a trabajar
% con el contexto an·logo a primer orden.

\section{Introduction}
Shelah and Stern proved that First Order Logic is not a good framework to study (from a Model-Theoretic point of view) classes of Metric Structures -such as Banach spaces-. In fact, they proved that these classes have a behavior similar to Second Order Logic with Predicades (see~\cite{ShSt}). Because of that, there was necessary to study a logic which would satisfy good properties as in First Order Logic (e.g., Downward %\linebreak
L\"owenheim-Skolem-Tarski Theorem, Compactness Theorem, etc.) and at the same time that would be a suitable logic to study classes of Metric Structures from this point of view. This was the beginning of {\it Continuous Logic} (for short, CL; see~\cite{CoMon}).
\\ \\
Shelah's Presentation Theorem is a very interesting and important
result in (discrete) {it Abstract Elementary Classes} (for short, AECs; a good framework to study classes of discrete structures which are not axiomatizable in First Order Logic) because this result allows us to work with Ehrenfeucht-Mostowski
models, and also it allows to prove the existence of arbitrarily large-enough
models in $\C{K}$ via the existence of Hanf numbers (see \cite{BaMon}). In fact, this follows from the fact that AECs are Projective Classes (the statement of the classical Shelah's Presentation Theorem in discrete AECs) and that AECs are controlled in some way by an infinitary logic. 
\\ \\
%\indent In that case, under categoricity, such models are used for proving stability in smaller cardinalities and also -without categoricity assumptions- those models allow us to calculate theHanf number of $\C{K}$ .
%\\ \\
{\it Metric Abstract Elementary Classes} (for short, MAECs) is a generalization of
{\it Continuous Logic} parallel to the notion of {\it Abstract Elementary Class}
 (see \cite{BaMon}), but we work with the completion of the
union of a elementary chain instead of working just with such union,
and also we work in that setting with density character instead of
cardinality. We follow the definitions and terminology given by \r{A}sa
Hirvonen and Tapani Hyttinen (\cite{Hi}). 
%
%But we have to point out that we do not consider 
%uniformly continuous functions, in our analysis it is enough to consider
%closed functions.
\\ \\
Hirvonen and Hyttinen (see \cite{Hi}) proved the following weaker
version of Shelah's Presentation Theorem

\begin{theorem}
Let $(\C{K},\ordencea)$ an MAEC of $L$-structures with
$|L|+LS^d(\C{K})\le \aleph_0$. Then for each $M\in \C{K}$ we can
define an expansion $M^*$ with Skolem functions $F_n^k$
($k,n<\omega$) such that:
\begin{enumerate}
\item If $A\subset M^*$ and $A$ is closed under the functions
$F_n^k$ then $\overline{A}\upharpoonright L\in \C{K}$ and
$\overline{A}\upharpoonright L\ordencea M$.
\item For all $a\in M$, $A_a:=\{(F^{length(a)}_n)^{M^*}(a):n<\omega
\}$ is such that
\begin{enumerate}
    \item $\overline{A_a}\upharpoonright L\in \C{K}$ and
    $\overline{A_a}\upharpoonright L\ordencea M$,
    \item If $b\subseteq a$ (as sets) then $b\in A_b \subset A_a$.
\end{enumerate}
\end{enumerate}
\end{theorem}

%Nuevo: enero 26 de 2009
However, they do not prove that an MAEC is a Projective Class with
omitting types. They used their version of Shelah's Presentation
Theorem for constructing Ehrenfeucht-Mostowski models in this
setting.
\\ \\
In \cite{Za} we refined their argument, providing an explicit
theory in Continuous Logic and an explicit set of types which work in a similar way as
in the original proof in the AEC setting, proving that an MAEC is in
fact a Projective Class. In \cite{Za} we claimed that the interpretations of the given function symbols in the extended language were not necessarily
uniformly continuous, which is a requirement for fitting this theory and this set of omitted types in the setting of Continuous Logic.
\\ \\
In this paper, we will prove that those interpretations are actually uniformly continuous (claim~\ref{UniformCont}).  For the sake of completeness, we will provide the proof given in \cite{Za}.
\\ \\
For basic notions of Continuous Logic, we refer the reader to\linebreak
\cite{CoMon}. For basic notions of MAEC, we refer the reader to
\cite{Hi,Za}.

\section{Metric Abstract Elementary Classes}
\begin{remark}
Through this paper, we call a {\it complete metric space} an
$L$-structure (in the context of {\it Continuous Logic}). %But we point out that we do not consider uniformly continuous symbols. For our purposes,  it is enough to consider closed functions. Uniformly continuity is used in (first order) Continuous Logic to prove the axiomatizability of  several classes of metric structures, but we are not interested to study this kind of classes.
\end{remark}

\begin{definition}\label{DensityCharacter}
The density character of a topological space is the smallest
cardinality of a dense subset of the space. If $X$ is a topological
space, we denote its density character by $dc(X)$. If $A$ is a
subset of a topological space $X$, we define
$dc(A):=dc(\overline{A})$.
\end{definition}

We consider a natural adaptation of the notion of {\it Abstract
Elementary Class} (see \cite{Gr} and \cite{BaMon}), but working in
the context of Continuous Logic (see \cite{CoMon}). We follow the
definitions given by \r{A}sa Hirvonen and Tapani Hyttinen (see \cite
{Hi}).

%Añadido agosto 6 de 2011
\begin{definition}\label{Structures}
Let $L$ be a language as in \cite{CoMon}, but without the uniform continuity modulus. A multi-sorted {\it metric $L$-structure} is a tuple
\[M:=(\{\langle A_i,d_i\rangle\}_{i\in I}, \mathbb{R}, \{c_j\}_{j\in J},\{F_k\}_{k\in K},\{R_l\}_{l\in L} ) \text{, where:}\]
\begin{enumerate}
\item Each $(A_i,d_i)$ is a complete metric space.
\item $\mathbb{R}$ is an isomorphic copy of the real field $(\mathbb{R}, +,\cdot, 0, 1, \le)$.
\item Each $c_j$ is a constant in a fixed sort $A_{i(j)}$.
\item Each $R_k$ is a continuous predicate; i.e., $R_k$ corresponds to a function $R_k:A_{i(k)_1}\times \cdots \times A_{i(k)_n}\to [0,1]$ which is closed, i.e.: if $(\overline{x})_{n<\omega}\to \overline{x}$ as tuples, then $(R_k(\overline{x}_n))_{n<\omega}\to R_k(\overline{x})$,  where $n$ is called the {\it arity of $R_k$}.
\item Each $F_l$ is  a function $F_l:A_{i(l)_1}\times \cdots \times A_{i(l)_m}\to A_{i(l)}$ which is closed %\marginpar{?`SER\'A EL MEJOR NOMBRE?}
(i.e.: if $(\overline{x})_{n<\omega}\to \overline{x}$ as tuples, then $(F_l(\overline{x}_n))_{n<\omega}\to F_l(\overline{x})$), where $m$ is called the {\it arity of $F_l$} .
\end{enumerate}

If it is clear that we are working in a metric context, we just called them {\it $L$-structures}.

\end{definition}

\begin{definition}\label{MAEC}
Let $\mathcal{K}$ be a class of $L$-structures as defined in \ref{Structures} above and $\ordencea$ be a binary relation defined in
$\C{K}$. We say that $(\C{K},\ordencea)$ is a {\it Metric Abstract
Elementary Class} (shortly {\it MAEC}) if:

\begin{enumerate}
\item $\C{K}$ and $\ordencea$ are closed under $\cong$.
\item $\ordencea$ is a partial order in $\C{K}$.
\item If $\mathfrak{M}\ordencea \mathfrak{N}$ then $\mathfrak{M}\subseteq \mathfrak{N}$.
\item (Tarski-Vaught chains) If $(\mathfrak{M}_i:i<\lambda)$ is a $\ordencea$-increasing chain then
    \begin{enumerate}
    \item the function symbols in $L$ can be uniquely interpreted on the completion of
        $\bigcup_{i<\lambda} \mathfrak{M}_i$ such that
        $\overline{\bigcup_{i<\lambda}\mathfrak{M}_i} \in \C{K}$
    \item for each $j < \lambda$ , $\mathfrak{M}_j \ordencea
        \overline{\bigcup_{i<\lambda} \mathfrak{M}_i}$
    \item if each $\mathfrak{M}_i\in \C{K}\in \mathfrak{N}$, then
        $\overline{\bigcup_{i<\lambda} \mathfrak{M}_i} \ordencea \mathfrak{N}$.
    \end{enumerate}
    \item (coherence) if $\mathfrak{M}_1\subseteq \mathfrak{M}_2\ordencea \mathfrak{M}_3$ and
        $\mathfrak{M}_1\ordencea \mathfrak{M}_3$, then $\mathfrak{M}_1\ordencea\mathfrak{M}_2$.
    \item (Downward L\"owenheim-Skolem) There exists a cardinality $LS^d(K)$ (which is called {\it L\"o\-wen\-heim-Skolem number}) such that if $\mathfrak{M} \in \C{K}$
        and $A \subseteq M$, then there exists $\mathfrak{N}\in \C{K}$ such that $dc(\mathfrak{N}) \le dc(A) + LS^d(K)$
        and $A\subseteq \mathfrak{N} \ordencea \mathfrak{M}$.
\end{enumerate}

\end{definition}

\begin{examples}
\begin{enumerate}
\item Any Continuous Logic Elementary Class with the elementary substructure
relation is an MAEC.
\item Positive bounded theories, where $\ordencea$ is interpreted by the approximate
elementary submodel relation (see \cite{HeIo}).
\item Compact Abstract Theories, see \cite{Be1,Be2}
\item The class of Banach spaces, where $\ordencea$ is interpreted
by the closed subspace relation (see \cite{Hi}).
\item Hilbert Spaces expanded with an unbounded closed selfadjoint operator (see~\cite{Ar}). This example is not necessarily axiomatizable in Continuous Logic.
\end{enumerate}
\end{examples}

\subsection{Some basic facts}
In this section, we mention some basic (and classic) facts towards getting a proof of Shelah Presentation Theorem. This basic facts are also used in the classic proof in the (discrete) Abstract Elementary Classes (for short, AECs), but for the sake of completeness we provide their statements.

\begin{fact}\label{DirOrder}
Let $(I,\le)$ be a directed partial order of size $\lambda$. Then
there exists a family $\{I_\alpha : \alpha<\lambda\}$ of suborders
of $I$ such that:
\begin{enumerate}
    \item Each $I_\alpha$ is a directed order and $|I_\alpha|<\lambda$
    \item If $\alpha<\beta<\lambda$, then $I_\alpha\le I_\beta$
    \item $I=\bigcup_{\alpha<\lambda} I_\alpha$.
\end{enumerate}
\end{fact}
\bsdem \cite{Ma}\edem

We prove the following fact in a similar way as in (discrete) AECs (mutatis mutandis). In fact,  we strongly use the Tarski-Vaught chains axiom (MAEC axiom). Notice that in MAECs, this axiom involves not just the union of the $\ordencea$-chain, we have to take the completion of that union.  Despite of the sketch of the proof is almost the same as in (dicrete) AECs, for the sake of completeness we provide a proof of this fact.

\begin{proposition}\label{DirSystemMAECs}
Let $(I,\le )$ be a directed partial order and $(\mathfrak{M}_i :i\in I)$ a $\ordencea$-directed system. Then:
\begin{enumerate}
    \item[(a)] $\overline{\bigcup_{i\in I}\mathfrak{M}_i} \in \C{K}$.
    \item[(b)] $\mathfrak{M}_j\ordencea \overline{\bigcup_{i\in I}\mathfrak{M}_i}$ for each $j\in I$.
    \item[(c)] If $\mathfrak{N}\in \C{K}$ and $\mathfrak{M}_j\ordencea \mathfrak{N}$ for each $j\in I$,
        then $\overline{\bigcup_{i\in I}\mathfrak{M}_i}\ordencea \mathfrak{N}$.
\end{enumerate}
\end{proposition}

\bdem 
\ Assume this fact holds for $\alpha<|I|$. By fact \ref{DirOrder} we have that there exists a
family $\{I_\alpha : \alpha<\lambda\}$ of suborders of $I$ such that:
\begin{enumerate}
    \item Each $I_\alpha$ is a directed order and $|I_\alpha|<|I|$
    \item If $\alpha<\beta<|I|$, then $I_\alpha\le I_\beta$
    \item $I=\bigcup_{\alpha<|I|} I_\alpha$.
\end{enumerate}

Define $\mathfrak{M}_{\alpha}:=\overline{\bigcup_{i\in I_\alpha} \mathfrak{M}_i}$. By induction hypothesis (b)
we have that $\mathfrak{M}_j\ordencea \mathfrak{M}_\alpha$ for every $j\in I_\alpha$. If $\alpha<\beta$,
since $I_\alpha\subseteq I_\beta$ then $\mathfrak{M}_j\ordencea \mathfrak{M}_\beta$ for every $j\in I_\alpha$.
By induction hypothesis (c) we have that $\mathfrak{M}_\alpha:=\overline{\bigcup_{j\in I_{\alpha}}
\mathfrak{M}_j}\ordencea \mathfrak{M}_\beta$.
\\\\
It is easy to check that $\overline{\bigcup_{\alpha<|I|}\mathfrak{M}_\alpha} =
\overline{\bigcup_{i\in I}\mathfrak{M}_i}$, so by definition \label{MAEC} 4) (a) we have that
$\overline{\bigcup_{i\in I}\mathfrak{M}_i} = \overline{\bigcup_{\alpha<|I|}\mathfrak{M}_\alpha} \in \C{K}$. Then
(a) holds.

If $j\in I$, there exists $\alpha<|I|$ such that $j\in I_\alpha$, so $\mathfrak{M}_j\ordencea \mathfrak{M}_\alpha$
(by induction hypothesis) and by definition \ref{MAEC} (4) (b) $\mathfrak{M}_\alpha\ordencea
\overline{\bigcup_{\alpha<|I|}\mathfrak{M}_\alpha} =
\overline{\bigcup_{i\in I}\mathfrak{M}_i}$.
Therefore $\mathfrak{M}_j\ordencea \overline{\bigcup_{i\in I} \mathfrak{M}_{i}}$, i.e. (b) holds.

Let $\mathfrak{N}$ be an $L$-structure in $\C{K}$ such that $\mathfrak{M}_j\ordencea \mathfrak{N}$ for each
$j\in I$. By induction hypothesis (c), for each $\alpha<|I|$we have that
$\mathfrak{M}_{\alpha}\ordencea \mathfrak{N}$. So, by definition \ref{MAEC} (4) (c) we have
that $\overline{\bigcup_{i\in I}\mathfrak{M}_i}=\overline{\bigcup_{\alpha<|I|}\mathfrak{M}_\alpha}
\ordencea \mathfrak{N}$. So, (c) holds.
\edem[Proposition~\ref{DirSystemMAECs}]

\begin{definition}[directed system]\label{Directed_System}
Let $\C{K}$ be a Category. A functor $D:(I,\le)\to \C{C}$ is said to be a {\it directed system} if and only if $(I,\le)$ is a directed ordered set. Set $M_k:=D(k)$ for every $k\in I$ and $f_{i,j}:=D((i,j)):M_i\to M_j$ the morphism associate to the unique $I$-morphism $(i,j):i\to j$ via $D$ whenever $i\le_I j$.
\end{definition}

\begin{definition}[directed limits]
We say that $\C{K}$ is {\it closed under directed limits} iff
for every directed system $D:(I,\le)\to \C{C}$ there exist $M\in ob(\C{K})$ and $\C{C}$-morphisms $f_{i,\infty}: M_i\to M$ ($i\in I$) such that 
\begin{enumerate}
\item for any $i\le_I j$ we have $f_{i,\infty}=f_{j,\infty}\circ f_{i,j}$
\item if any $N\in ob(\C{C})$ has a system of $\C{C}$-morphisms $g_{i,\infty}:M_i\to N$ which safisties 1. above, then there exists a unique $\C{C}$-morphism $h:M\to N$ such that $g_{i,\infty}=h\circ f_{i,\infty}$.
\[
\begin{diagram}
\node{M_i} \arrow{e,t}{f_{i,\infty}} \arrow{se,b}{g_{i,\infty}}
\node{M}\arrow{s,r}{h}\\
\node[2]{N}
\end{diagram}
\]
\end{enumerate}
Such morphisms $f_{i,\infty}$ are called {\it canonical morphisms}.
\end{definition}

\begin{corollary}\label{Dir_Sys}
An MAEC $\C{K}$ (viewed as a category with morphisms the\linebreak
$\C{K}$-embeddings) is closed under directed limits.
\end{corollary}
\bdem
 \ This is a direct consequence of proposition~\ref{DirSystemMAECs} and axiom 1. of definition~\ref{MAEC}.
\edem[Corollary~\ref{Dir_Sys}]

\section{The main question.}
%\begin{question}
%Does the Shelah presentation theorem hold for an MAEC?
%\end{question}
%\marginpar{ESCRIBIR el teorema y preliminares en AEC}
%We will intend to adapt the proof of Shelah's Presentation Theorem which holds in the context of
%{\it Abstract Elementary Classes}, but we will indicate the problems which we will have when
%we do this.

%AÒadido diciembre 28 de 2008
\begin{definition}
Let $\C{K}$ be a class of $L$-models in the continuous logic
setting (but with closed functions instead of uniformly continuous functions). We say that $\C{K}$ is a {\it projective class with omitting types} (shortly, PC$\Gamma$ class) iff there
exist an expansion $L'$ of $L$, an $L'$-theory $T'$ and a set of
$L'$-types $\Gamma$ such that
$\C{K}=PC_L(T',\Gamma):=\{\mathfrak{M}\upharpoonright L :
\M{M}\models T' \text{\ and } \M{M} \text{\ omits every type in }
\Gamma\}$.
\end{definition}

And finally, we provide a proof of the version of Shelah Presentation Theorem in the setting of MAECs. We have to clarify that although the sketch of the proof is almost the same of the discrete proof, we are working in a metric setting and we have to change lots of details in the proof. This is the proof given in \cite{Za}, which we provide for the sake of completeness. In \cite{Za} we claimed that the interpretations of the function symbols of the Skolemization are not necessarily uniformly continuous, but in this paper we will prove that they are in fact uniformly continuous (claim~\ref{UniformCont}). 

\begin{theorem}[Shelah's Presentation Theorem in MAECs]\label{PresentationTh}
Given an MAEC $(\C{K},\ordencea)$, there exist an expansion $L'$ of
$L(\C{K})$, a $L'$-theory $T'$ and a set of $T'$-types $\Gamma$ such
that $\mathcal{K}=PC_L(T',\Gamma)$ (i.e.: $\C{K}$ is a projective
class with omitting types).
\end{theorem}

\bdem
Let $L'$ be the language obtained from $L(\mathcal{K})$ by adding
new $n$-ary function symbols $F^n_i$ ($i<LS(\mathcal{K})$). Let $T'$
be the theory which says that $\sup x_0 \dots \sup x_{n-1}
|F^n_i(x_o,\cdots,x_{n-1})-x_i|=0$ for all $i<n$.  Notice that all $F_i^n$ are defined as projections, so they are continuous (and therefore, closed).
%$\forall x_0 \dots \forall x_{n-1} F^n_i(x_o,\cdots,x_{n-1})=x_i$ for all $i<n$

Take $\mathfrak{M}\models T'$ and $\overline{a}\in M$. For $\overline{b}<_{subtuple}\overline{a}$,
define $U_{\overline{b}}:=\{F^m_i(\overline{b}) : i<LS(\mathcal{K})\}$, where $m:=l(\overline{b})$.
\\\\
Define $\Gamma$ as follows: for each $\mathfrak{M}'\models T'$ and each tuple
$\overline{a}\in M'$, $tp(\overline{a}/\emptyset,\mathfrak{M}')\in \Gamma$ {\bf unless}
we have the following two conditions:

\begin{enumerate}
    \item Taking $\overline{b}\le_{subtuple} \overline{a}$, $U_{\overline{a}}$
        is a dense subset of a submodel of $\mathfrak{M}'\upharpoonright L$,
        which we denote by $\mathfrak{M}_{\overline{b}}$, and $\mathfrak{M}_{\overline{b}}\in \mathcal{K}$.
    \item Taking $\overline{b}\le_{subtuple} \overline{a}$, we have
        $\mathfrak{M}_{\overline{b}}\ordencea \mathfrak{M}_{\overline{a}}$.
\end{enumerate}

Take $\mathfrak{M}\in PC_L(T',\Gamma)$. So, there exists
$\mathfrak{M}'\models T'$ such that omits all the types in $\Gamma$
and $\mathfrak{M}=\mathfrak{M}'\upharpoonright L$. Consider the sets
$U_{\overline{a}}$, for each $\overline{a}\in M$. As $\mathfrak{M}'$
omits all the types in $\Gamma$, each $U_{\overline{a}}$ is the
universe of a submodel $\mathfrak{M}_{\overline{a}}$ of
$\mathfrak{M}$ such that $\mathfrak{M}_{\overline{a}}\in
\mathcal{K}$. By proposition \ref{DirSystemMAECs}
%\marginpar{ARREGLAR ESA PRUEBA ... VALEN TODAS LAS CONDICIONES, NO LA ⁄NICA QUE SE PUSO},
we have that $\overline{\bigcup_{\overline{a}\in M} \mathfrak{M}_{\overline{a}}}\in \mathcal{K}$.

Since $M'\models T'$, we have that $\overline{a}\in U_{\overline{a}}$ and so
$\overline{\bigcup_{\overline{a}\in M} \mathfrak{M}_{\overline{a}}}=\mathfrak{M}$. So, $\mathfrak{M}\in \mathcal{K}$.
\\ \\
In the other way, take $\mathfrak{M}\in \mathcal{K}$. We define $\mathfrak{M}'$ as follows:
for $n=0$, choose $\mathfrak{M}_{\emptyset}\ordencea \mathfrak{M}$ of density character
$LS^d(\mathcal{K})$ and let $U_{\emptyset}:=\{ F^0_i : i<LS^d(\mathcal{K}) \}$ be an enumeration of
a dense subset of $M_{\emptyset}$.
Having done this for $n$, let $\overline{a}\in M$ be of lenght $n+1$. Choose
$\mathfrak{M}_{\overline{a}}\ordencea \mathfrak{M}$ of density character $LS^d(\mathcal{K})$ which contains
$\bigcup\{M_{\overline{b}}:\overline{b}<_{subtuple}\overline{a}\}$ and $\overline{a}$, and let
$U_{\overline{a}}:=\{F^{n+1}_i(\overline{a}):i<LS(\mathcal{K})\}$ be an enumeration of
a dense subset of $M_{\overline{a}}$ such that $F^{n+1}_i(\overline{a})=a_i$ for $i<n+1$, where $\overline{a}:=(a_0,\cdots,a_{n})$.

\begin{claim}\label{UniformCont}
$F^{n+1}_i$ is an uniformly continuous function.
\end{claim}
\bdem
\ Notice that $F^{n+1}_i:\mathfrak{M}^{n+1}\to \mathfrak{M}$ defined as above is a projection. The topology in $\mathfrak{M}^{n+1}$ is given by the metric defined by $d^*(\overline{b},\overline{c}):=\max_{j<n+1} \{ d(b_j,c_j) \}$ (where $\overline{b}:=(b_0,\cdots, b_n)$, $\overline{c}:=(c_0,\cdots,c_n)$ and $d$ is the metric of $\mathfrak{M}$). Let $\E>0$. Taking $\delta:=\epsilon$, if $d^*(\overline{b},\overline{c}):=\max_{j<n+1}\{d(b_j,c_j)\}<\delta$, then
\begin{eqnarray*} 
d(F^{n+1}_i(\overline{b}),F^{n+1}_i(\overline{c})) &=& d(b_i,c_i)\\
&\le& \max_{j<n+1}\{d(b_j,c_j)\}\\
&=:& d^*(\overline{b},\overline{c})\\
&<& \delta\\
&:=& \E
\end{eqnarray*}
\edem[Claim~\ref{UniformCont}]

Notice that $\mathfrak{M}'\models T'$ and $\mathfrak{M}'$ omits every type
in $\Gamma$.
%\\ \\
\begin{center}
\scalebox{0.5} % Change this value to rescale the drawing.
{
\begin{pspicture}(0,-4.4735594)(15.289075,4.4535594)
\definecolor{color55}{rgb}{0.8,0.0,0.0}
\definecolor{color56}{rgb}{0.0,0.6,0.2}
\psframe[linewidth=0.04,dimen=outer](15.22,4.4535594)(0.0,-2.1664405)
\psline[linewidth=0.04cm](0.06,-4.3664403)(15.14,-4.2864404)
\scalebox{3}{\psdots[dotsize=0.12,dotstyle=diamond*](2.33,-1.45)}
\psbezier[linewidth=0.04,linecolor=color55](0.02,1.253)(0.02,0.45356)(7.52,3.18)(8.52,3.15356)(9.52,3.13)
(15.27,2.3452)(15.14,1.35356)
\scalebox{2}{\rput(0.4,0.98){$F_i^n$}}
\psline[linewidth=0.04cm,linecolor=color56](7,-4.3064404)(7,4.4135594)
\scalebox{2}{\rput(4,0.02){$\mathfrak{M}_{\overline{a}}$}}
\scalebox{2}{\rput[c](3.5,-2.5){$\overline{a}$}}
\scalebox{2}{\rput[c](3,-3){$n$-tuples in $\mathfrak{M}$}}
\scalebox{2}{\rput(7.86,1.24){$\mathfrak{M}$}}
\scalebox{3}{\psdots[dotsize=0.12](2.15,0.95)}
\scalebox{2}{\rput(2.5,1.6){$F^{n}_i(\overline{a})$}}
\end{pspicture}
}
\end{center}
\ 
\edem[Theorem~\ref{PresentationTh}]
%\ \\ \\
%\begin{question}[Berenstein]
%Is $F^{n+1}_i$ an uniformly continuous function? Actually, it does
%not matter if this question is not true because in Continuous Logic
%we need uniform continuity for proving the {\it Compactness
%Theorem}, but in a general way we do not need that requirement. In our setting, we just require that the function symbols are weakly closed.
%\end{question}

%We weak a little the Berenstein's questions, and we can assure that in fact these functions are continuous because all projections are continuous.
%\\ \\
%\indent Uniform continuity of function and relational symbols is used for
%proving a version of compactness theorem in Continuous Logic, see
%\cite{CoMon}. But we do not require that conditions because we are
%interested to work in a general setting where even that theorem does
%not hold. For all our further studies -see for example~\cite{ViZa15}- we just require that every symbol function is closed.

%Pregunta numero de Hanf
\begin{remark}
Let $L$ be a (first order) language and $\mu:=\sup\{|L|,\kappa\}$. Since Shelah proved that any AEC $\C{K}$ is a $PC\Gamma(\kappa,2^\kappa)$ class (i.e., $\C{K}=PC_L(T,\Gamma 
)$ with $|T|=\kappa$ and $|\Gamma|=2^\kappa$) with $LS(\C{K})=\kappa$, as a consequence of the M. Morley's omitting types theorem -see \cite{Morley}- we have that there exists a cardinality $H_1:=\beth_{(2^{\mu})^+}$ such that if $\C{K}$ is an AEC of $L$-structures with $LS(\C{K})=\kappa$ %\warning
%$PC\Gamma(\kappa,2^\kappa)$ class of $L$-structures  
such that if there exists $M\in \C{K}$ of cardinality $>H_1$ then there exists a model in $\C{K}$ in any cardinality $>H_1$. %As a consequence of Morley's result, we have that PC classes have a Hanf number as well.
%Since Shelah proved that AECs are $PC\Gamma(\kappa,2^\kappa)$ classes, then each AEC has a Hanf number. 
But its existence strongly depends on the existence of Hanf numbers in (discrete) omiting type classes (which depends on infinitary logics, see \cite{Morley}). Despite there are some recent works about some intends of providing suitable notions of metric infinitary logics (see~\cite{Ea}), it is still open if this can yields a suitable analysis of Hanf numbers which implies that metric PC classes have a Hanf number, and so MAECs does as well. However, W. Boney proved that Hanf numbers exist for MAECs by using a kind of adjoint functors between the original MAEC and an auxiliary discrete AEC defined by some well-behaved dense subsets closed under function symbols (see~\cite{Bo}).
\end{remark}

%\\ \\
%\indent Xavier Caicedo talked (in the Universidad Nacional deColombia Model Theory Seminar in 2009) about a very interesting relation between uniform continuity and compactness. He showed that if we take a connective $\square_f$ associate to funcion $f:[0,1]^n\to [0,1]$, then the extension $\L (\square_f)$ of the \L ukasiewicz logic  (see \cite{Ha}) adding that connective $\square_f$ is compact if and only if $f$ is continuous.


\begin{thebibliography}{9}
\bibitem[Ar14]{Ar} C. Argoty, {\it Model theory of a Hilbert space expanded with an unbounded closed selfadjoint operator}, Math. Logic Q. 60 (6) pp 403--424, 2014. 
\bibitem[Ba09]{BaMon} J. Baldwin, \textit{Categoricity in Abstract Elementary Classes}, monograph, http://www.math.uic.edu/\verb"~"jbaldwin/pub/AEClec.pdf
\bibitem[Be03]{Be1} I. Ben-Yaacov, \textit{Positive Model Theory and compact abstract
theories}, J. Math. Log. 3 (1) pp. 85-118, 2003.
\bibitem[Be05]{Be2} I. Ben-Yaacov, \textit{Uncountable dense categoricity in
cats}, J. Symbolic Logic  70 (3), pp. 829-860, 2005.
\bibitem[BeBeHeUs08]{CoMon} I. Ben-Yaacov, A. Berenstein, C. W. Henson, A. Usvyatsov,
    {\it Model theory for metric structures}, Model theory with Applications to
    Algebra and Analysis Vol. 2 (Eds. Z. Chatzidakis, D. Macpherson,
    A. Pillay, A. Wilkie) London Math Soc. Lecture Note Series 350, Cambridge Univ Press, 2008.
\bibitem[Bo14]{Bo} W. Boney, {\it A presentation theorem for continuous logic and Metric Abstract Elementary Classes}, preprint.  arXiv:1408.3624 
\bibitem[Ea14]{Ea} C. Eagle, {\it Omitting types for infinitary $[0,1]$-valued logic}, Ann. Pure Appl. Logic 165 pp. 913--932, 2014. 
\bibitem[Gr02]{Gr} R. Grossberg, {\it Classification theory for abstract elementary classes},
    Logic and Algebra, ed. Yi Zhang, Contemporary Mathematics 302, AMS pp. 165--204, 2002.
\bibitem[HeIo02]{HeIo} C.W. Henson, J. Iovino, {\it Ultraproducts in
analysis}, in Analysis and logics (mons, 1997), London
Math. Soc. Lecture Note Ser.  262 pp. 1--110, Cambridge University Press,
2002.
\bibitem[Ha98]{Ha} P. H\'ajek, {\it Mathematics of fuzzy logics}, Kluwer, 1998
\bibitem[HiHy08]{Hi} \r{A}. Hirvonen, T. Hyttinen, {\it Categoricity in homogeneous complete metric spaces}, Arch. Math. Logic 48 pp. 269--322, 2009. 
\bibitem[Ma85]{Ma} J.A. Makowsky, {\it Abstract Embedding Relations}, in Model - theoretic
logics, ed. by J. Barwise and S. Feferman, Springer - Verlag,
New York, 1985.
\bibitem[Mo65]{Morley} M.~Morley, {\it Omiting classes of elements}, in {\em The theory of models}, Henkin, Addison and Tarski (ed), pp. 265--273. North Holland, 1965.
\bibitem[ShSt78]{ShSt} S. Shelah, J. Stern, {\it The Hanf number of the First Order  Theory of Banach Spaces}, Trans. American Math. Soc. 244 pp 147--171, 1978.
\bibitem[Za11]{Za} P. Zambrano, {Around Superstability in Metric Abstract Elementary Classes}. Ph.D. Thesis, Universidad Nacional de Colombia, 2011.
\end{thebibliography}
\end{document}